\documentclass{article}
\usepackage{mitpress}
\usepackage{graphicx}
\usepackage{amsmath}
\usepackage{amsfonts}
\usepackage{amssymb}
\newtheorem{theorem}{Theorem}

\newdimen\dummy
\dummy=\oddsidemargin
\graphicspath{ {Result_astro_128/}  {Result_astro_256/}  {Result_edges/} {Result_Kon_astro_(chains)/} {Result_Lena/}}

\begin{document}

\title{Polyharmonic Daubechies type wavelets in Image Processing and Astronomy, II}
\author{Ognyan Kounchev, Damyan Kalaglarsky, Milcho Tsvetkov}
\maketitle

\textbf{Abstract: }\emph{We consider the application of the polyharmonic
subdivision wavelets to Image Processing, in particular to Astronomical
Images. The results show an essential advantage over some standard
multivariate wavelets and a potential for better compression.}


\section{Introduction}

In \cite{KounchevKalaglarsky2010} we have constructed the basic elements of
the polyharmonic Wavelet Analysis. There we have provided all details of the
construction of the filters for the whole infinite family of father and mother
wavelets. Here we will finalize the multivariate construction and provide
experimental evidence of the power of the new wavelets.

\section{Polyharmonic Wavelet Analysis \label{SpolyharmonicSubdivision}}

The \emph{Polyharmonic Wavelet Analysis} arises in the context of the
polyharmonic subdivision. In general, it is obtained when the
(one-dimensional) polynomials of degree $2N-1$ are replaced by polyharmonic
functions of degree $N.$ The main element of the construction is the
application of the Fourier transform in $n$ variables in $\mathbb{R}^{n+1}.$
For simplicity we assume that the Image which we consider is a function
$u\left(  t,y\right)  ,$ for $t\in\mathbb{R}$ and $2\pi-$periodic in
$y\in\mathbb{R}^{n}.$ We have
\begin{align}
u\left(  t.y\right)   &  =\sum_{\eta\in\mathbb{Z}^{n}}\widehat{u}_{\eta
}\left(  t\right)  e^{i\left\langle \eta,y\right\rangle }\qquad\text{where the
Fourier coefficients are}\label{FourierTransform}\\
\widehat{u}_{\eta}\left(  t\right)   &  =\frac{1}{\left(  2\pi\right)  ^{n}%
}\int_{0}^{2\pi}\cdot\cdot\cdot\int_{0}^{2\pi}u\left(  t.y\right)
e^{-i\left\langle \eta,y\right\rangle }dy.\nonumber
\end{align}
For every fixed $\eta\in\mathbb{Z}^{n}$ and $\xi=\left|  \eta\right|  $ we
make Wavelet Analysis which is based on the non-stationary father and mother
wavelets $\left\{  \varphi_{m}^{\xi}\left(  t\right)  \right\}  _{m\geq0}$ and
$\left\{  \psi_{m}^{\xi}\left(  t\right)  \right\}  _{m\in\mathbb{Z}}$ of
Theorem 4. Recall that they have refinement coefficients $M_{j}^{\left[
m\right]  ,\xi}$ and $\left(  -1\right)  ^{j}M_{1-j}^{\left[  m\right]  ,\xi}$
respectively, where the polynomial $M^{\left[  m\right]  ,\xi}\left(
z\right)  $ was defined at the end of \cite{KounchevKalaglarsky2010}. This
gives the wavelet expansion
\begin{equation}
u\left(  t.y\right)  =\sum_{m,j\in\mathbb{Z}}\sum_{\eta\in\mathbb{Z}^{n}%
}\gamma_{\eta,m,j}\psi_{m}^{\left|  \eta\right|  }\left(  t-j\right)
e^{i\left\langle \eta,y\right\rangle }. \label{waveletExpansion}%
\end{equation}

We summarize the results in the following theorem, for the proofs we refer to
\cite{DynKounchevLevinRender}.

\begin{theorem}
For every $\xi\geq0$ the functions $\psi_{m}^{\xi}\left(  t\right)  $ generate
a non-stationary MRA of $L_{2}\left(  \mathbb{R}\right)  $ by means of the
definitions
\begin{align*}
V_{m}^{\xi}  &  :=clos_{L_{2}\left(  \mathbb{R}\right)  }\left\{  \varphi
_{m}^{\xi}\left(  t-j\right)  :j\in\mathbb{Z}\right\} \\
L_{2}\left(  \mathbb{R}\right)   &  =clos_{L_{2}\left(  \mathbb{R}\right)
}\left(  \bigcup_{m\in\mathbb{Z}}V_{m}^{\left|  \xi\right|  }\right)  .
\end{align*}
Here $clos_{L_{2}\left(  \mathbb{R}\right)  }$ denotes the usual closure in a
linear and topological sense with respect to the space $L_{2}\left(
\mathbb{R}\right)  .$ Respectively, the Wavelet Analysis spaces $W_{m}%
^{\left|  \xi\right|  }$ defined by means of
\[
W_{m}^{\xi}\bigoplus V_{m}^{\xi}=V_{m+1}^{\xi}%
\]
are generated as
\[
W_{m}^{\xi}=clos_{L_{2}\left(  \mathbb{R}\right)  }\left\{  \psi_{m}^{\xi
}\left(  t-j\right)  :j\in\mathbb{Z}\right\}  .
\]
\end{theorem}

As in the polyspline Wavelet Analysis studied in \cite{okbook}, the
one-dimensional MRA $V_{j}^{\left|  \eta\right|  }$ which is generated for
every $\eta\in\mathbb{Z}^{n}$ is related to a ''polyharmonic MRA'' $PV_{j}$ as
$PV_{j}=\bigoplus_{\eta\in\mathbb{Z}^{n}}V_{j}^{\left|  \eta\right|  }$ by
means of the Fourier transform formula (\ref{FourierTransform}). Respectively,
we have also the polyharmonic Wavelet Analysis defined by $PV_{j}\bigoplus
PW_{j}=PV_{j+1}$ and reduced to infinitely many one-dimensional Wavelet
Analyses for every $\eta\in\mathbb{Z}^{n}$ in the form $PW_{j}=\bigoplus
_{\eta\in\mathbb{Z}^{n}}W_{j}^{\left|  \eta\right|  }$ by means of the Fourier
transform formula (\ref{FourierTransform}).

Respectively, the space
\[
L_{2,per}\left(  \mathbb{R}^{n+1}\right)  =\left\{  f\left(  t,y\right)
:t\in\mathbb{R},\ y\in\mathbb{R}^{n};f\text{ is }2\pi-\text{periodic in
}y\right\}
\]
is generated by an infinite amount of \textbf{orthogonal (non-stationary)
mother wavelets}
\begin{equation}
\left\{  \psi_{m}^{\left|  \eta\right|  }\left(  t\right)  e^{i\left\langle
\eta,y\right\rangle }:\eta\in\mathbb{Z}^{n},\ m\in\mathbb{Z}\right\}  .
\label{wavelets}%
\end{equation}
Although these wavelets do not have a compact support in direction $y,$ it is
important to note that they have ''elongated support'' in the $y$ direction,
which plays an essential role for the sparseness of the ''edge
representation'' in the images, a point of view nicely described in
\cite{Donoho}, see also \cite{mallat}.

\section{The algorithm}

From the above it follows that the algorithm for Decomposition and
Reconstruction of the \textbf{Polyharmonic Subdivision Wavelet Analysis
(PhSdWav)} consists of a Fourier transform in direction $y$ and consecutive
application, for every $\eta\in\mathbb{Z}^{n},$ of one-dimensional Wavelet
Analysis in direction $t.$ It consists of the following steps:

\begin{enumerate}
\item  As usually we take as \textbf{first approximation} to the image the
image itself where the expansion is in terms of the scaling functions
$\left\{  \varphi_{m}^{\left|  \eta\right|  }\left(  t\right)  \right\}
_{m\geq0}.$

\item  Transforming to Wavelet Domain - this is the main difference with the
usual approach. Roughly speaking, we obtain $\widehat{u}_{\eta}\left(
t\right)  $ in formula (\ref{FourierTransform}) and make its Wavelet Analysis
for every $\eta\in\mathbb{Z}^{n}$ which gives us the coefficients
$\gamma_{\eta,m,j}$ in the wavelet expansion (\ref{waveletExpansion}). Let us
note that these coefficients are in general complex numbers.

\item  Filtering the small (insignificant) coefficients by using a threshold parameter.

\item  Quantizing the remaining coefficients.

\item  Encoding by arithmetic coder.
\end{enumerate}

The algorithm was written in Matlab.

\section{Experiments with images}

We apply our Wavelet Analysis to different classes of images -- synthetic and
real. We compare the results of our method  PhSdWav with parameter $N=9$ (see
\cite{KounchevKalaglarsky2010}) against the $2d$ wavelet transform with
Daubechies wavelets $db9$ (which are tensor product wavelets and correspond to
the one-dimensional Daubechies wavelets with parameter $N=9$). We use the
standard function $dwt2$ in Matlab which contains a large class of $2d$
wavelets, and $db9$ is one of the most used representatives. \textbf{In all
Figures below} we first display on top-left the original image, on top-right
is the result of PhSdWav compressed representation obtained by wavelet
thresholding, and on the bottom is the result of the $db9$ compressed
representation obtained by  wavelet thresholding. In order to make comparison
between our method PhSdWav and $db9$ we perform  thresholding which results in
equal PSNR. 

\subsection{Elongated support tests on $64\times64$ pixel edges\label{Sedges}}

For the proper compression and representation of images it is of fundamental
importance to understand how the method applies to ''edges''. For that reason
we will first show the performance of our method to some synthetic edge images.

\subsubsection{Vertical edge}

In this test the image is a \textbf{vertical edge} which is a function
$f_{vert}\left(  t,y\right)  $ defined in the rectangle $\left[  0,1\right]
\times\left[  0,1\right]  $ by
\[
f_{vert}\left(  t,y\right)  :=\left\{
\begin{array}
[c]{c}%
0\qquad\text{for }0\leq t\leq1/2\\
1\qquad\text{for }1/2<t\leq1
\end{array}
\right.  .
\]
We take a sampling of $f\left(  t,y\right)  $ by a $64\times64$ pixel image.
The performance of the polyharmonic Wavelets and the $db9$ wavelets is
provided in \textbf{Table 1} and \textbf{Figure 1}  below.%

\begin{align*}
&  \qquad\qquad\qquad\qquad\qquad\qquad\qquad\qquad\qquad\qquad\text{Table
1.}\\
&
\begin{tabular}
[c]{||l||l||l||l||}\hline\hline
& Num. coeffs & PSNR & Entropy\\\hline\hline
Pol. Wav. & \textbf{48} & \textbf{120.7379} & \textbf{0}\\\hline\hline
db9 & \textbf{960} & \textbf{121.7412} & \textbf{1.01e+003}\\\hline\hline
\end{tabular}
\end{align*}

\begin{figure}[th]
\centering
\includegraphics[width=0.80\textwidth
]{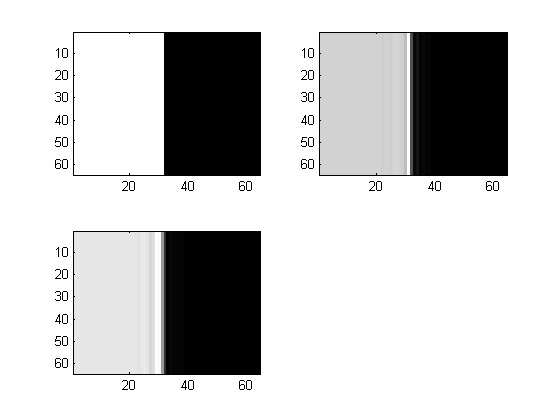}%
\caption{Vertical Edge}%
\label{fig:Result_vertical}%
\end{figure}

As we said above the wavelets of our basis (\ref{wavelets}) show elongation in
direction $y$ and this makes the above image $f_{vert}$ perfect for
representation and compression with the polyharmonic subdivision wavelets. We
see that the quality of representation by means of PhSdWav (measured by the
PSNR) is the same as the quality of representation by means of the standard
$2-$dimensional Daubechies wavelet $db9$ but the number of coefficients is
$20$ times less!

The image $f_{vert}$ would correspond to the one-dimensional image provided by
the Heaviside function, $f\left(  t\right)  =0$ for $0\leq t\leq1/2$ ,
$f\left(  t\right)  =1$ for $1/2<t\leq1.$ The wavelets have an excellent
performance for such images.

Just the opposite is the next example.

\subsubsection{Horizontal edge}

In this test the image is a \textbf{horizontal edge}, which is a function
$f_{hor}\left(  t,y\right)  $ defined in the rectangle $\left[  0,1\right]
\times\left[  0,1\right]  $ by
\[
f_{hor}\left(  t,y\right)  :=\left\{
\begin{array}
[c]{c}%
0\qquad\text{for }0\leq y\leq1/2\\
1\qquad\text{for }1/2<y\leq1
\end{array}
\right.  .
\]
We take a sampling of $f_{hor}\left(  t,y\right)  $ by a $64\times64$ pixel
image. The performance of the polyharmonic Wavelets and the $db9$ wavelets is
provided in \textbf{Table 2} and \textbf{Figure 2} below.%

\begin{align*}
&  \qquad\qquad\qquad\qquad\qquad\qquad\qquad\qquad\qquad\qquad\text{Table
2.}\\
&
\begin{tabular}
[c]{||l||l||l||l||}\hline\hline
& Num. coeffs & PSNR & Entropy\\\hline\hline
Pol. Wav. & \textbf{2560} & \textbf{121.4740} & \textbf{2.52e+003}%
\\\hline\hline
db9 & \textbf{960} & \textbf{121.7412} & \textbf{1.01e+003}\\\hline\hline
\end{tabular}
\end{align*}

\begin{figure}[th]
\centering
\includegraphics[width=0.80\textwidth
]{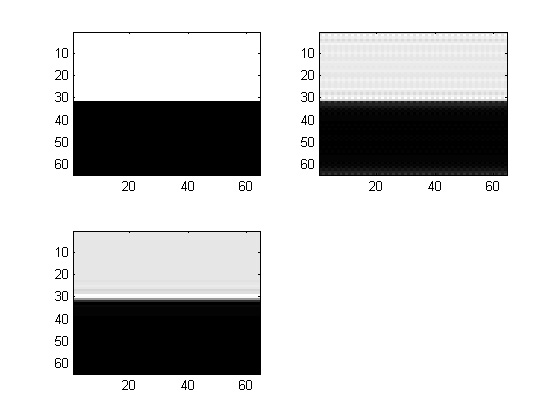}%
\caption{Horizontal Edge}%
\label{fig:Result_horizontal}%
\end{figure}

We see now how inefficient the PhSdWav are for the representation of $f_{hor}$
compared to $db9$ wavelets, since for the same PSNR we need almost $3$ times
more coefficients. This result is expected since PhSdWav has a very strong
vertical orientation.

\subsubsection{Skewed edge}

Somewhere in the middle is the case of a skewed edge. In this test the image
of \textbf{skewed edge}, which is a function $f_{skew}\left(  t,y\right)  $
defined in the rectangle $\left[  0,1\right]  \times\left[  0,1\right]  $ by
(this is an approximate definition)
\[
f_{skew}\left(  t,y\right)  :=\left\{
\begin{array}
[c]{c}%
0\qquad\text{for }0\leq t\leq\frac{1}{6}\left(  y+2.5\right)  \\
1\qquad\text{for }\frac{1}{6}\left(  y+2.5\right)  <t\leq1
\end{array}
\right.  .
\]
We take a sampling of $f_{skew}\left(  t,y\right)  $ by a $64\times64$ pixel
image. The performance of the polyharmonic Wavelets and the $db9$ wavelets is
provided in \textbf{Table 3} and \textbf{Figure 3 }below.%

\begin{align*}
&  \qquad\qquad\qquad\qquad\qquad\qquad\qquad\quad\quad\qquad\qquad\text{Table
3.}\\
&
\begin{tabular}
[c]{||l||l||l||l||}\hline\hline
& Num. coeffs & PSNR & Entropy\\\hline\hline
Pol. Wav. & \textbf{512} & \textbf{126.8204} & \textbf{ 386.7158}%
\\\hline\hline
db9 & \textbf{940} & \textbf{121.9320} & \textbf{977.3171}\\\hline\hline
\end{tabular}
\end{align*}

\begin{figure}[th]
\centering
\includegraphics[width=0.80\textwidth
]{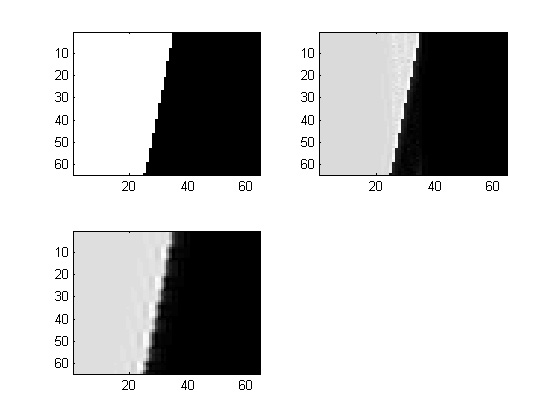}\caption{Skewed
Edge}%
\label{fig:Result_skewed}%
\end{figure}

We see that although the edge is prettily skew (the slope is approximately
$\tan\alpha=6$ ) the PhSdWav have a very nice performance compared to $db9.$

\subsection{Lena image}

Let us make a more detailed analysis of the seminal image of Lena in
$128\times128$ pixels. We provide the Compression Ratio and the Bits Encoded
for both methods in \textbf{Table 4} and \textbf{Figure 4 }below.
\begin{align*}
&  \qquad\qquad\qquad\qquad\qquad\qquad\quad\qquad\qquad\qquad\qquad
\qquad\qquad\qquad\qquad\text{Table 4.}\\
&
\begin{tabular}
[c]{||l||l||l||l||l||l||}\hline\hline
& N. coeffs & PSNR & Entropy & Cmpr. rat. & Bits Enc.\\\hline\hline
Pol. Wav. & \textbf{4,150} & \textbf{77.7338} & \textbf{5.15e+004} &
\textbf{5.9815} & \textbf{42,802}\\\hline\hline
db9 & \textbf{5,187} & \textbf{77.6742} & \textbf{5.65e+004} & \textbf{5.9562}%
& \textbf{44,012}\\\hline\hline
\end{tabular}
\end{align*}

\begin{figure}[th]
\centering
\includegraphics[width=1.00\textwidth
]{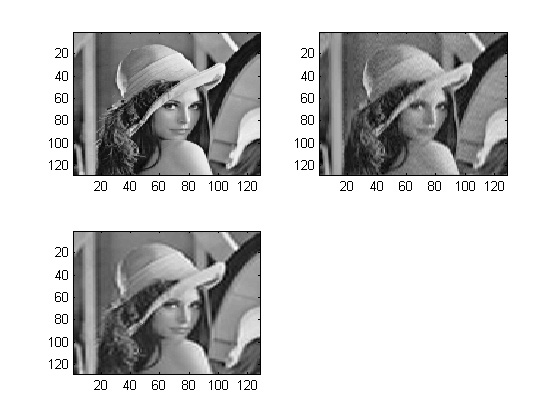}\caption{Lena Image}%
\label{fig:Result_Lena}%
\end{figure}

We see that for the same PSNR our method PhSdWav is better -- it needs $20\%$
less coefficients. This is an astonishingly good result since as we have seen
in the experiments with the synthetic edges in Section \ref{Sedges}, our
method of PhSdWav is adapted to the vertical direction $y$ and the Lena image
does not have a dominating orientation of the edges.

\subsection{Astronomical Images}

Hereafter we analyze different types of astronomical images.

\subsubsection{The Pleiades star cluster, a part with faint objects only,
sized $128\times128$}

First we consider an image which is a part of the Pleiades image containing
only faint stars. The results are provided in \textbf{Table 5} and
\textbf{Figure 5 }below.%
\begin{align*}
&  \qquad\qquad\qquad\qquad\qquad\qquad\qquad\qquad\qquad\qquad\qquad
\qquad\qquad\qquad\text{Table 5.}\\
&
\begin{tabular}
[c]{||l||l||l||l||l||l||}\hline\hline
& N. coeffs & PSNR & Entropy & Cmpr. rat. & Bits enc.\\\hline\hline
Pol. Wav. & \textbf{2,663} & \textbf{77.5805} & \textbf{2.93e+004} &
\textbf{7.6028} & \textbf{33,456}\\\hline\hline
db9 & \textbf{5,187} & \textbf{76.3970} & \textbf{4.52e+004} & \textbf{6.9317}%
& \textbf{37,818}\\\hline\hline
\end{tabular}
\end{align*}

\begin{figure}[h]
\centering
\includegraphics[width=1.00\textwidth
]{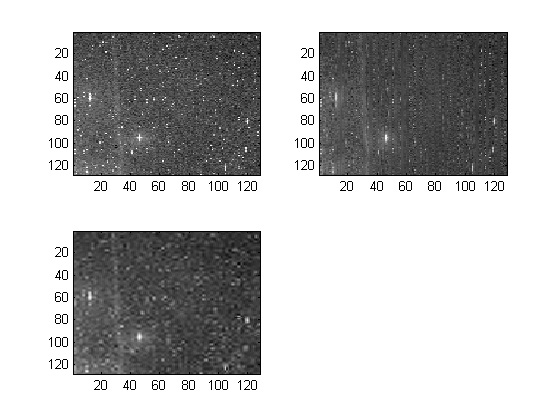}%
\caption{Pleiades star cluster, faint objects, 128x128 pixels}%
\label{fig:Result_astro_128}%
\end{figure}

We see that for (almost) the same PSNR we have much less coefficients, the
compression ratio is higher, and our entropy is much less, which shows a
better potential for compression.

\subsubsection{The Pleiades star cluster, a part with bright objects, sized
$256\times256$}

Next we consider a part of the Pleiades containing several bright objects. The
results are provided in \textbf{Table 6 }and \textbf{Figure 6} below%
\begin{align*}
&  \qquad\qquad\qquad\qquad\qquad\qquad\qquad\qquad\qquad\qquad\qquad
\qquad\qquad\qquad\text{Table 6.}\\
&
\begin{tabular}
[c]{||l||l||l||l||l||l||}\hline\hline
& N. coeffs & PSNR & Entropy & Cmpr. rat. & Bits enc.\\\hline\hline
Pol. Wav. & \textbf{15,159} & \textbf{50.0717} & \textbf{2.89e+005} &
\textbf{13.8338} & \textbf{73,750}\\\hline\hline
db9 & \textbf{18,502} & \textbf{49.4214} & \textbf{3.57e+005} &
\textbf{8.5379} & \textbf{122,814}\\\hline\hline
\end{tabular}
\end{align*}

\begin{figure}[h]
\centering
\includegraphics[width=1.00\textwidth
]{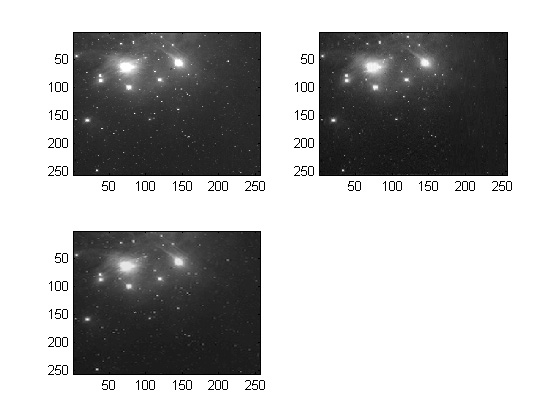}%
\caption{Pleiades star cluster, bright objects, 256x256 pixels}%
\label{fig:Result_astro_256}%
\end{figure}

Again we see that for the same PSNR we have not only less coefficients but
also a much higher compression ratio.

\subsubsection{The Pleiades star cluster, chain, sized $128\times128$}

Here we consider a chain image of a part of the Pleiades image. The results
are provided in \textbf{Table 7} and \textbf{Figure 7 }below%
\begin{align*}
&  \qquad\qquad\qquad\qquad\qquad\qquad\qquad\qquad\qquad\qquad\quad
\qquad\qquad\qquad\qquad\text{Table 7.}\\
&
\begin{tabular}
[c]{||l||l||l||l||l||l||}\hline\hline
& N. coeffs & PSNR & Entropy & Cmpr. rat. & Bits enc.\\\hline\hline
Pol. Wav. & \textbf{2,764} & \textbf{88.6818} & \textbf{2.4849e+004} &
\textbf{8.8706} & \textbf{28,528}\\\hline\hline
db9 & \textbf{5,340} & \textbf{87.5730} & \textbf{5.7468e+004} &
\textbf{8.4636} & \textbf{30,973}\\\hline\hline
\end{tabular}
\end{align*}

\begin{figure}[h]
\centering
\includegraphics[width=1.00\textwidth
]{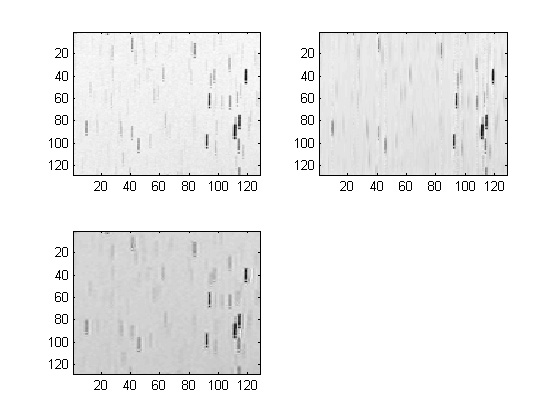}%
\caption{Pleiades star cluster, chain image, 128x128 pixels}%
\label{fig:Result}%
\end{figure}

For the same PSNR, the number of coefficients of our method PhSdWav is twice
less than the number of  $db9$ coefficients. This result is expected since the
image is vertically oriented. However the compression ratio is not much bigger
apparently due to more noise in the background of the image.

\textbf{Acknowledgement.} The first named author was sponsored partially by
the Alexander von Humboldt Foundation, and all authors were sponsored by
Project DO--2-275/2008 ''Astroinformatics'' with Bulgarian NSF.

\textbf{ABOUT THE AUTHORS}

Ognyan Kounchev, Prof. Dr., Institute of Mathematics and Informatics,
Bulgarian Academy of Science, tel. +359-2-9793851; kounchev@gmx.de

Damyan Kalaglarsky, Institute of Astronomy, Bulgarian Academy of Science, tel.
+359-2-9793851; damyan@skyarchive.org.

Milcho Tsvetkov, Assoc. Prof., Dr., Institute of Astronomy, Bulgarian Academy
of Science, tel. +359-2-9795935; milcho@skyarchive.org.
\end{document}